\documentclass[12pt,vatola]{article}

\usepackage{graphicx}

\def\c{\centerline}

\def\re#1{\par\hangindent\parindent\indent\llap{#1\enspace}\ignorespaces}

\def\no{\noindent}

\begin{document}

\c{\bf\large Pseudo-Manifold Geometries with Applications}

\vskip 6mm

\c{Linfan Mao}\vskip 2mm

\c{\scriptsize (Chinese Academy of Mathematics and System Science,
Beijing 100080, P.R.China)}

\c{\scriptsize E-mail: maolinfan@163.com}

\vskip 6mm

\begin{minipage}{130mm}

\no{\bf Abstract}: {\small  A Smarandache geometry is a geometry
which has at least one Smarandachely denied axiom($1969$), i.e.,
an axiom behaves in at least two different ways within the same
space, i.e., validated and invalided, or only invalided but in
multiple distinct ways and a Smarandache $n$-manifold is a
$n$-manifold that support a Smarandache geometry. Iseri provided a
construction for Smarandache $2$-manifolds by equilateral
triangular disks on a plane and a more general way for Smarandache
$2$-manifolds on surfaces, called {\it map geometries} was
presented by the author in $[9]-[10]$ and $[12]$. However, few
observations for cases of $n\geq 3$ are found on the journals. As
a kind of Smarandache geometries, a general way for constructing
dimensional $n$ pseudo-manifolds are presented for any integer
$n\geq 2$ in this paper. Connection and principal fiber bundles
are also defined on these manifolds. Following these
constructions, nearly all existent geometries, such as those of
{\it Euclid geometry, Lobachevshy-Bolyai geometry, Riemann
geometry, Weyl geometry, K\"{a}hler geometry} and {\it Finsler
geometry}, ...,etc., are their sub-geometries. }\vskip 2mm

\no{\bf Key Words}: {\small Smarandache geometry, Smarandache
manifold, pseudo-manifold, pseudo-manifold geometry,
multi-manifold geometry, connection, curvature, Finsler geometry,
Riemann geometry, Weyl geometry and K\"{a}hler geometry.}\vskip
2mm

\no{\bf AMS(2000)}: {\small 51M15, 53B15, 53B40, 57N16}

\end{minipage}

\vskip 8mm

\no{\bf \S $1.$ Introduction}

\vskip 4mm

\no Various geometries are encountered in update mathematics, such
as those of {\it Euclid geometry, Lobachevshy-Bolyai geometry,
Riemann geometry, Weyl geometry, K\"{a}hler geometry} and {\it
Finsler geometry}, ..., etc.. As a branch of geometry, each of
them has been a kind of spacetimes in physics once and contributes
successively to increase human's cognitive ability on the natural
world. Motivated by a combinatorial notion for sciences: {\it
combining different fields into a unifying field}, Smarandache
introduced {\it neutrosophy and neutrosophic logic} in references
$[14]-[15]$ and Smarandache geometries in $[16]$.

\vskip 4mm

\no{\bf Definition $1.1$}([8][16]) \ {\it An axiom is said to be
Smarandachely denied if the axiom behaves in at least two
different ways within the same space, i.e., validated and
invalided, or only invalided but in multiple distinct ways.

A Smarandache geometry is a geometry which has at least one
Smarandachely denied axiom($1969$).}

\vskip 3mm

\no{\bf Definition $1.2$} \ {\it For an integer $n, n\geq 2$, a
Smarandache $n$-manifold is a $n$-manifold that support a
Smarandache geometry.}

\vskip 3mm

Smarandache geometries were applied to construct many world from
conservation laws as a mathematical tool([2]). For Smarandache
$n$-manifolds, Iseri constructed Smarandache manifolds for $n=2$
by equilateral triangular disks on a plane in $[6]$ and $[7]$ (see
also $[11]$ in details). For generalizing Iseri's Smarandache
manifolds, {\it map geometries} were introduced in $[9]-[10]$ and
$[12]$, particularly in $[12]$ convinced us that these map
geometries are really Smarandache $2$-manifolds. Kuciuk and
Antholy gave a popular and easily understanding example on an
Euclid plane in $[8]$. Notice that in $[13]$, these {\it
multi-metric space} were defined, which can be also seen as
Smarandache geometries. However, few observations for cases of
$n\geq 3$ and their relations with existent manifolds in
differential geometry are found on the journals. The main purpose
of this paper is to give general ways for constructing dimensional
$n$ pseudo-manifolds for any integer $n\geq 2$. Differential
structure, connection and principal fiber bundles are also
introduced on these manifolds. Following these constructions,
nearly all existent geometries, such as those of {\it Euclid
geometry, Lobachevshy-Bolyai geometry, Riemann geometry, Weyl
geometry, K\"{a}hler geometry} and {\it Finsler geometry},
...,etc., are their sub-geometries.

Terminology and notations are standard used in this paper. Other
terminology and notations not defined here can be found in these
references $[1],[3]-[5]$.

For any integer $n, n\geq 1$, an {\it $n$-manifold} is a Hausdorff
space $M^n$, i.e., a space that satisfies the $T_2$ separation
axiom, such that for $\forall p\in M^n$, there is an open
neighborhood $U_p, p\in U_p\subset M^n$ and a homeomorphism
$\varphi_p: U_p\rightarrow {\bf R}^n$ or ${\bf C}^n$,
respectively.

Considering the differentiability of the homeomorphism $\varphi:
U\rightarrow {\bf R}^n$ enables us to get the conception of
differential manifolds, introduced in the following.

An {\it differential $n$-manifold} $(M^n, {\mathcal A})$ is an
$n$-manifold $M^n, M^n=\bigcup\limits_{i\in I}U_i$, endowed with a
$C^r$ differential structure ${\mathcal
A}=\{(U_{\alpha},\varphi_{\alpha})| \alpha\in I\}$ on $M^n$ for an
integer $r$ with following conditions hold.

$(1)$ \ $\{U_{\alpha}; \alpha\in I\}$ is an open covering of
$M^n$;

$(2)$ \ For $\forall \alpha,\beta\in I$, atlases
$(U_{\alpha},\varphi_{\alpha})$ and $(U_{\beta},\varphi_{\beta})$
are {\it equivalent}, i.e., $U_{\alpha}\bigcap
U_{\beta}=\emptyset$ or $U_{\alpha}\bigcap
U_{\beta}\not=\emptyset$ but the {\it overlap maps}

$$
\varphi_{\alpha}\varphi_{\beta}^{-1}:
\varphi_{\beta}(U_{\alpha\bigcap U_{\beta}})\rightarrow
\varphi_{\beta}(U_{\beta}) \ \ {\rm and} \ \
\varphi_{\beta}\varphi_{\alpha}^{-1}:
\varphi_{\beta}(U_{\alpha\bigcap U_{\beta}})\rightarrow
\varphi_{\alpha}(U_{\alpha})
$$

\no are $C^r$;

$(3)$ \ ${\mathcal A}$ is maximal, i.e., if $(U,\varphi)$ is an
atlas of $M^n$ equivalent with one atlas in ${\mathcal A}$, then
$(U,\varphi)\in{\mathcal A}$.

An $n$-manifold is {\it smooth} if it is endowed with a
$C^{\infty}$ differential structure. It is well-known that a
complex manifold $M^n_c$ is equal to a smooth real manifold
$M^{2n}_r$ with a natural base

$$\{\frac{\partial}{\partial x^i}, \frac{\partial}{\partial
y^i}| \ 1\leq i\leq n\}$$

\no for $T_pM^n_c$, where $T_pM^n_c$ denotes the tangent vector
space of $M^n_c$ at each point $p\in M^n_c$.

\vskip 8mm

\no{\bf \S $2.$ Pseudo-Manifolds}

\vskip 4mm

\no These Smarandache manifolds are non-homogenous spaces, i.e.,
there are singular or inflection points in these spaces and hence
can be used to characterize warped spaces in physics. A
generalization of ideas in map geometries can be applied for
constructing dimensional $n$ pseudo-manifolds.

\vskip 4mm

\no{\bf Construction $2.1$} \ {\it Let $M^n$ be an $n$-manifold
with an atlas ${\mathcal A}=\{(U_p, \varphi_p)| p\in M^n\}$. For
$\forall p\in M^n$ with a local coordinates
$(x_1,x_2,\cdots,x_n)$, define a spatially directional mapping
$\omega: p\rightarrow {\bf R}^n$ action on $\varphi_p$ by

$$\omega: p\rightarrow\varphi_p^{\omega}(p)=\omega(\varphi_p(p))=
(\omega_1,\omega_2,\cdots,\omega_n),$$

\no i.e., if a line $L$ passes through $\varphi(p)$ with direction
angles $\theta_1,\theta_2,\cdots,\theta_n$ with axes ${\bf
e}_1,{\bf e}_2,\cdots,{\bf e}_n$ in ${\bf R}^n$, then its
direction becomes

$$\theta_1-\frac{\vartheta_1}{2}+\sigma_1,\theta_2-\frac{\vartheta_2}{2}+\sigma_2,
\cdots, \theta_n-\frac{\vartheta_n}{2}+\sigma_n$$

\no after passing through $\varphi_p(p)$, where for any integer
$1\leq i\leq n$, $\omega_i\equiv\vartheta_i(mod 4\pi)$,
$\vartheta_i\geq 0$ and

\[
\sigma_i=\left\{\begin{array}{cc}
\pi,&  if \ \ 0\leq\omega_i< 2\pi,\\
0, & if \ \ 2\pi < \omega_i< 4\pi.
\end{array}
\right.
\]

\no A manifold $M^n$ endowed with such a spatially directional
mapping $\omega: M^n\rightarrow {\bf R}^n$ is called an
$n$-dimensional pseudo-manifold, denoted by $(M^n, {\mathcal
A}^{\omega})$.}

\vskip 4mm

\no{\bf Theorem $2.1$} \ {\it For a point $p\in M^n$ with local
chart $(U_p,\varphi_p)$, $\varphi_p^{\omega}=\varphi_p$ if and
only if $\omega(p)=(2\pi k_1,2\pi k_2,\cdots,2\pi k_n)$ with
$k_i\equiv 1(mod 2)$ for $1\leq i\leq n$.}

\vskip 3mm

{\it Proof} \  By definition, for any point $p\in M^n$, if
$\varphi_p^{\omega}(p)=\varphi_p(p)$, then
$\omega(\varphi_p(p))=\varphi_p(p)$. According to Construction
$2.1$, this can only happens while $\omega(p)=(2\pi k_1,2\pi
k_2,\cdots,$ $2\pi k_n)$ with $k_i\equiv 1(mod 2)$ for $1\leq
i\leq n$. \ \ $\natural$

\vskip 4mm

\no{\bf Definition $2.1$} \ {\it A spatially directional mapping
$\omega: M^n\rightarrow {\bf R}^n$ is euclidean if for any point
$p\in M^n$ with a local coordinates $(x_1,x_2,\cdots,x_n)$,
$\omega(p)=(2\pi k_1,2\pi k_2,\cdots,2\pi k_n)$ with $k_i\equiv
1(mod 2)$ for $1\leq i\leq n$, otherwise, non-euclidean.}

\vskip 3mm

\no{\bf Definition $2.2$} \ {\it Let $\omega: M^n\rightarrow {\bf
R}^n$ be a spatially directional mapping and $p\in (M^n, {\mathcal
A}^{\omega})$, $\omega(p)(mod
4\pi)=(\omega_1,\omega_2,\cdots,\omega_n)$. Call a point $p$
elliptic, euclidean or hyperbolic in direction ${\bf e}_i$, $1\leq
i\leq n$ if $o\leq\omega_i< 2\pi$, $\omega_i=2\pi$ or $2\pi <
\omega_i < 4\pi$. }

\vskip 3mm

Then we get a consequence by Theorem $2.1$.

\vskip 4mm

\no{\bf Corollary $2.1$} \ {\it Let $(M^n, {\mathcal A}^{\omega})$
be a pseudo-manifold. Then $\varphi_p^{\omega}=\varphi_p$ if and
only if every point in $M^n$ is euclidean.}

\vskip 4mm

\no{\bf Theorem $2.2$} \ {\it Let $(M^n, {\mathcal A}^{\omega})$
be an $n$-dimensional pseudo-manifold and $p\in M^n$. If there are
euclidean and non-euclidean points simultaneously or two elliptic
or hyperbolic points in a same direction in $(U_p, \varphi_p)$,
then $(M^n, {\mathcal A}^{\omega})$ is a Smarandache
$n$-manifold.}

\vskip 3mm

{\it Proof} \ On the first, we introduce a conception for locally
parallel lines in an $n$-manifold. Two lines $C_1, C_2$ are said
{\it locally parallel} in a neighborhood $(U_p,\varphi_p)$ of a
point $p\in M^n$ if $\varphi_p(C_1)$ and $\varphi_p(C_2)$ are
parallel straight lines in ${\bf R}^n$.

In $(M^n, {\mathcal A}^{\omega})$, the axiom that {\it there are
lines pass through a point locally parallel a given line} is
Smarandachely denied since it behaves in at least two different
ways, i.e., {\it one parallel, none parallel}, or {\it one
parallel, infinite parallels}, or {\it none parallel, infinite
parallels}.

If there are euclidean and non-euclidean points in $(U_p,
\varphi_p)$ simultaneously, not loss of generality, we assume that
$u$ is euclidean but $v$ non-euclidean, $\omega(v)(mod
4\pi)=(\omega_1,\omega_2,\cdots,\omega_n)$ and
$\omega_1\not=2\pi$. Now let $L$ be a straight line parallel the
axis ${\bf e}_1$ in ${\bf R}^n$. There is only one line $C_u$
locally parallel to $\varphi^{-1}_p(L)$ passing through the point
$u$ since there is only one line $\varphi_p(C_q)$ parallel to $L$
in ${\bf R}^n$ by these axioms for Euclid spaces. However, if $0<
\omega_1<2\pi$, then there are infinite many lines passing through
$u$ locally parallel to $\varphi_p^{-1}(L)$ in $(U_p,\varphi_p)$
since there are infinite many straight lines parallel $L$ in ${\bf
R}^n$, such as those shown in Fig.$2.1(a)$ in where each straight
line passing through the point $\overline{u}=\varphi_p(u)$ from
the shade field is parallel to $L$.

\includegraphics[bb=0 5 100 130]{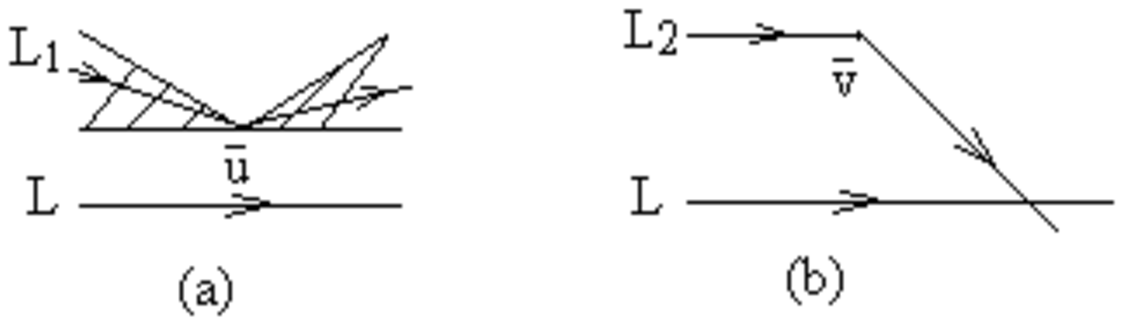}

\c{\bf Fig.$2.1$}\vskip 4mm

\no But if $2\pi < \omega_1< 4\pi$, then there are no lines
locally parallel to $\varphi_p^{-1}(L)$ in $(U_p,\varphi_p)$ since
there are no straight lines passing through the point
$\overline{v}=\varphi_p(v)$ parallel to $L$ in ${\bf R}^n$, such
as those shown in Fig.$2.1(b)$.

\includegraphics[bb=0 5 100
140]{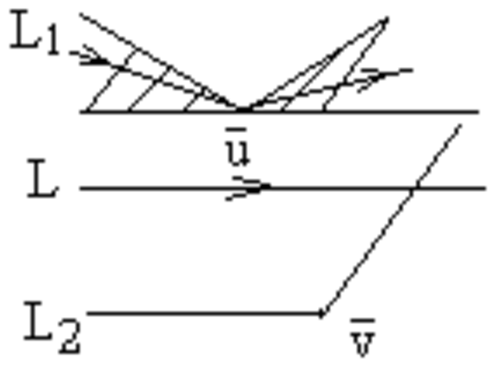}

\c{\bf Fig.$2.2$}

\vskip 2mm

If there are two elliptic points $u,v$ along a direction
$\overrightarrow{O}$, consider the plane ${\mathcal P}$ determined
by $\omega(u),\omega(v)$ with $\overrightarrow{O}$ in ${\bf R}^n$.
Let $L$ be a straight line intersecting with the line $uv$ in
${\mathcal P}$. Then there are infinite lines passing through $u$
locally parallel to $\varphi_p(L)$ but none line passing through
$v$ locally parallel to $\varphi_p^{-1}(L)$ in $(U_p,\varphi_p)$
since there are infinite many lines or none lines passing through
$\overline{u}=\omega(u)$ or $\overline{v}=\omega(v)$ parallel to
$L$ in ${\bf R}^n$, such as those shown in Fig.$2.2$.

Similarly, we can also get the conclusion for the case of
hyperbolic points. Since there exists a Smarandachely denied axiom
in $(M^n,{\mathcal A}^\omega)$, it is a Smarandache manifold. This
completes the proof. \ \ $\natural$

For an Euclid space ${\bf R}^n$, the homeomorphism $\varphi_p$ is
trivial for $\forall p\in {\bf R}^n$. In this case, we abbreviate
$({\bf R}^n,{\mathcal A}^\omega)$ to $({\bf R}^n, \omega)$.

\vskip 4mm

\no{\bf Corollary $2.2$} \ {\it For any integer $n\geq 2$, if
there are euclidean and non-euclidean points simultaneously or two
elliptic or hyperbolic points in a same direction in $({\bf R}^n,
\omega)$, then $({\bf R}^n, \omega)$ is an $n$-dimensional
Smarandache geometry.}

\vskip 3mm

Particularly, Corollary $2.2$ partially answers an open problem in
$[12]$ for establishing Smarandache geometries in ${\bf R}^3$.

\vskip 4mm

\no{\bf Corollary $2.3$} \ {\it If there are points $p, q\in {\bf
R}^3$ such that $\omega(p)(mod 4\pi)\not=(2\pi,2\pi,2\pi)$ but
$\omega(q)(mod 4\pi)=(2\pi k_1,2\pi k_2,2\pi k_3)$, where
$k_i\equiv 1(mod 2), 1\leq i\leq 3$ or $p,q$ are simultaneously
elliptic or hyperbolic in a same direction of ${\bf R}^3$, then
$({\bf R}^3, \omega)$ is a Smarandache space geometry.}

\vskip 3mm

\no{\bf Definition $2.3$} \ {\it For any integer $r\geq 1$, a
$C^r$ differential Smarandache $n$-manifold $(M^n, {\mathcal
A}^{\omega})$ is a Smarandache $n$-manifold $(M^n, {\mathcal
A}^{\omega})$ endowed with  a differential structure ${\mathcal
A}$ and a $C^r$ spatially directional mapping $\omega$. A
$C^{\infty}$ Smarandache $n$-manifold $(M^n, {\mathcal
A}^{\omega})$ is also said to be a smooth Smarandache
$n$-manifold.}

\vskip 3mm

According to Theorem $2.2$, we get the next result by definitions.

\vskip 4mm

\no{\bf Theorem $2.3$} \ {\it Let $(M^n,{\mathcal A})$ be a
manifold and $\omega: M^n\rightarrow {\bf R}^n$ a spatially
directional mapping action on ${\mathcal A}$. Then $(M^n,
{\mathcal A}^{\omega})$ is a $C^r$ differential Smarandache
$n$-manifold for an integer $r\geq 1$ if the following conditions
hold:}

$(1)$ \ {\it there is a $C^r$ differential structure ${\mathcal
A}=\{(U_{\alpha},\varphi_{\alpha})| \alpha\in I\}$ on $M^n$;}

$(2)$ \ {\it $\omega$ is $C^r$;}

$(3)$ \ {\it there are euclidean and non-euclidean points
simultaneously or two elliptic or hyperbolic points in a same
direction in $(U_p, \varphi_p)$ for a point $p\in M^n$.}

\vskip 3mm

{\it Proof} \  The condition $(1)$ implies that $(M^n, {\mathcal
A})$ is a $C^r$ differential $n$-manifold and conditions $(2)$,
$(3)$ ensure $(M^n, {\mathcal A}^{\omega})$ is a differential
Smarandache manifold by definitions and Theorem $2.2$. \ \
$\natural$

For a smooth differential Smarandache $n$-manifold $(M^n,
{\mathcal A}^{\omega})$, a function $f: M^n\rightarrow {\bf R}$ is
said smooth if for $\forall p\in M^n$ with an chart
$(U_p,\varphi_p)$,

$$f\circ(\varphi_p^{\omega})^{-1}: (\varphi_p^{\omega})(U_p)\rightarrow {\bf R}^n$$

\no is smooth. Denote by $\Im_p$ all these $C^{\infty}$ functions
at a point $p\in M^n$.

\vskip 4mm

\no{\bf Definition $2.4$} \ {\it Let $(M^n, {\mathcal
A}^{\omega})$ be a smooth differential Smarandache $n$-manifold
and $p\in M^n$. A tangent vector $v$ at $p$ is a mapping $v:
\Im_p\rightarrow {\bf R}$} with these following conditions hold.

$(1)$ \ $\forall g,h\in\Im_p, \forall\lambda\in {\bf R}, \
v(h+\lambda h)=v(g)+\lambda v(h);$

$(2)$ \ $\forall g,h\in\Im_p, v(gh)=v(g)h(p)+g(p)v(h).$

\vskip 3mm

Denote all tangent vectors at a point $p\in(M^n, {\mathcal
A}^{\omega})$ by $T_pM^n$ and define addition¡°+¡±and scalar
multiplication¡°$\cdot$¡±for $\forall u,v\in T_pM^n,
\lambda\in{\bf R}$ and $f\in\Im_p$ by

$$(u+v)(f)=u(f)+v(f), \ \ (\lambda u)(f)=\lambda\cdot u(f).$$

Then it can be shown immediately that $T_pM^n$ is a vector space
under these two operations¡°+¡±and¡°$\cdot$¡±.

Let $p\in(M^n, {\mathcal A}^{\omega})$ and $\gamma: (-\varepsilon,
\varepsilon)\rightarrow {\bf R}^n$ be a smooth curve in ${\bf
R}^n$ with $\gamma(0)=p$. In $(M^n, {\mathcal A}^{\omega})$, there
are four possible cases for tangent lines on $\gamma$ at the point
$p$, such as those shown in Fig.$2.3$, in where these bold lines
represent tangent lines.

\includegraphics[bb=0 5 100 120]{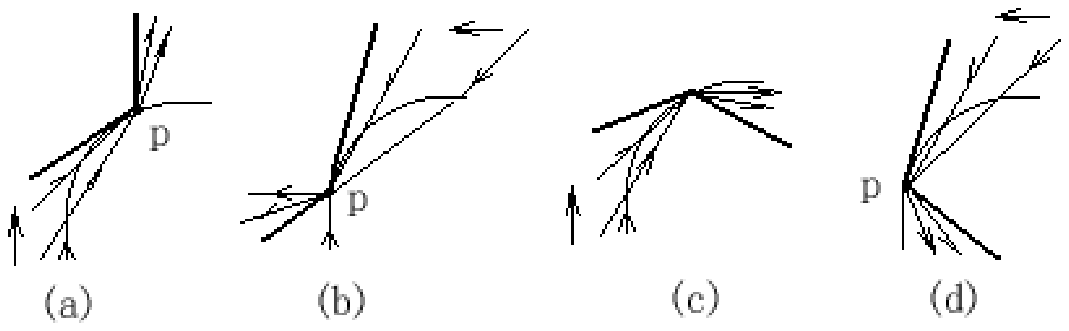}

\c{\bf Fig.$2.3$}

\vskip 2mm

By these positions of tangent lines at a point $p$ on $\gamma$, we
conclude that {\it there is one tangent line at a point $p$ on a
smooth curve if and only if $p$ is euclidean in $(M^n, {\mathcal
A}^{\omega})$}. This result enables us to get the dimensional
number of a tangent vector space $T_pM^n$ at a point $p\in (M^n,
{\mathcal A}^{\omega})$.

\vskip 4mm

\no{\bf Theorem $2.4$} \ {\it For any point $p\in(M^n, {\mathcal
A}^{\omega})$ with a local chart $(U_p, \varphi_p)$,
$\varphi_p(p)=(x_1^,x_2^0,\cdots,x_n^0)$, if there are just $s$
euclidean directions along ${\bf e}_{i_1},{\bf
e}_{i_2},\cdots,{\bf e}_{i_s}$ for a point , then the dimension of
$T_pM^n$ is}

$${\rm dim}T_pM^n \ = \ 2n-s$$

\no{\it with a basis}

$$\{\frac{\partial}{\partial x^{i_j}}|_p \ | \ 1\leq j\leq s\}
\bigcup\{\frac{\partial^-}{\partial x^l}|_p,
\frac{\partial^+}{\partial x^l}|_p \ | \ 1\leq l\leq n \ and \
l\not=i_j, 1\leq j\leq s \}.$$

\vskip 3mm

{\it Proof} \ We only need to prove that

$$\{\frac{\partial}{\partial x^{i_j}}|_p \ | \ 1\leq j\leq s\}
\bigcup\{\frac{\partial^-}{\partial x^l},
\frac{\partial^+}{\partial x^l}|_p \ | \ 1\leq l\leq n \ and \
l\not=i_j, 1\leq j\leq s \} \ \ (2.1)$$

\no is a basis of $T_pM^n$. For $\forall f\in\Im_p$, since $f$ is
smooth, we know that

\begin{eqnarray*}
f(x)&=&
f(p)+\sum\limits_{i=1}^n(x_i-x_i^0)\frac{\partial^{\epsilon_i}
f}{\partial x_i}(p)\\
&+&
\sum\limits_{i,j=1}^n(x_i-x_i^0)(x_j-x_j^0)\frac{\partial^{\epsilon_i}f}{\partial
x_i}\frac{\partial^{\epsilon_j}f}{\partial x_j}+R_{i,j,\cdots,k}
\end{eqnarray*}

\no for $\forall x=(x_1,x_2,\cdots,x_n)\in\varphi_p(U_p)$ by the
Taylor formula in ${\bf R}^n$, where each term in
$R_{i,j,\cdots,k}$ contains
$(x_i-x_i^0)(x_j-x_j^0)\cdots(x_k-x_k^0)$, $\epsilon_l\in\{+,-\}$
for $1\leq l\leq n$ but $l\not=i_j$ for $1\leq j\leq s$ and
$\epsilon_l$ should be deleted for $l=i_j, 1\leq j\leq s$.

Now let $v\in T_pM^n$. By Definition $2.4(1)$, we get that

\begin{eqnarray*}
v(f(x))&=&
v(f(p))+v(\sum\limits_{i=1}^n(x_i-x_i^0)\frac{\partial^{\epsilon_i}
f}{\partial x_i}(p))\\
&+&
v(\sum\limits_{i,j=1}^n(x_i-x_i^0)(x_j-x_j^0)\frac{\partial^{\epsilon_i}f}{\partial
x_i}\frac{\partial^{\epsilon_j}f}{\partial x_j})+
v(R_{i,j,\cdots,k}).
\end{eqnarray*}

Application of the condition $(2)$ in Definition $2.4$ shows that

$$v(f(p))=0, \ \ \sum\limits_{i=1}^nv(x_i^0)\frac{\partial^{\epsilon_i}
f}{\partial x_i}(p)=0,$$

$$v(\sum\limits_{i,j=1}^n(x_i-x_i^0)(x_j-x_j^0)\frac{\partial^{\epsilon_i}f}{\partial
x_i}\frac{\partial^{\epsilon_j}f}{\partial x_j})=0$$

\no and

$$v(R_{i,j,\cdots,k})=0.$$

\no Whence, we get that

$$v(f(x))=\sum\limits_{i=1}^nv(x_i)\frac{\partial^{\epsilon_i}
f}{\partial
x_i}(p)=\sum\limits_{i=1}^nv(x_i)\frac{\partial^{\epsilon_i}
}{\partial x_i}|_p(f). \ \ (2.2)$$

The formula $(2.2)$ shows that any tangent vector $v$ in $T_pM^n$
can be spanned by elements in $(2.1)$.

All elements in $(2.1)$ are linearly independent. Otherwise, if
there are numbers $a^1,a^2,\cdots,a^{s}, a_1^+,a_1^-,
a_2^+,a_2^-,\cdots,a_{n-s}^+,a_{n-s}^-$ such that

$$\sum\limits_{j=1}^{s}a_{i_j}\frac{\partial}{\partial x_{i_j}}
+\sum\limits_{i\not=i_1,i_2,\cdots,i_s, 1\leq i\leq
n}a_i^{\epsilon_i}\frac{\partial^{\epsilon_i}}{\partial
x_i}|_p=0,$$

\no where $\epsilon_i\in\{+,-\}$, then we get that

$$a_{i_j}=(\sum\limits_{j=1}^{s}a_{i_j}\frac{\partial}{\partial x_{i_j}}
+\sum\limits_{i\not=i_1,i_2,\cdots,i_s, 1\leq i\leq
n}a_i^{\epsilon_i}\frac{\partial^{\epsilon_i}}{\partial
x_i})(x_{i_j})=0$$

\no for $1\leq j\leq s$ and

$$a_i^{\epsilon_i}=(\sum\limits_{j=1}^{s}a_{i_j}\frac{\partial}{\partial x_{i_j}}
+\sum\limits_{i\not=i_1,i_2,\cdots,i_s, 1\leq i\leq
n}a_i^{\epsilon_i}\frac{\partial^{\epsilon_i}}{\partial
x_i})(x_i)=0$$

\no for $i\not=i_1,i_2,\cdots,i_s, 1\leq i\leq n$. Therefore,
$(2.1)$ is a basis of the tangent vector space $T_pM^n$ at the
point $p\in (M^n, {\mathcal A}^{\omega})$. \ \ $\natural$

Notice that ${\rm dim}T_pM^n=n$ in Theorem $2.4$ if and only if
all these directions are euclidean along ${\bf e}_1,{\bf
e}_2,\cdots,{\bf e}_n$. We get a consequence by Theorem $2.4$.

\vskip 4mm

\no{\bf Corollary $2.4$}([4]-[5]) \ {\it Let $(M^n,{\mathcal A})$
be a smooth manifold and $p\in M^n$. Then}

$${\rm dim}T_pM^n=n$$

\no{\it with a basis}

$$\{\frac{\partial}{\partial x^i}|_p \ | \ 1\leq i\leq n\}.$$

\vskip 3mm

\no{\bf Definition $2.5$} \ {\it For $\forall p\in(M^n, {\mathcal
A}^{\omega})$, the dual space $T_p^*M^n$ is called a co-tangent
vector space at $p$.}

\vskip 4mm

\no{\bf Definition $2.6$} \ {\it For $f\in\Im_p, d\in T_p^*M^n$
and $v\in T_pM^n$, the action of $d$ on $f$, called a differential
operator $d: \Im_p\rightarrow {\bf R}$, is defined by}

$$df \ = \ v(f).$$

\vskip 3mm

Then we immediately get the following result.

\vskip 4mm

\no{\bf Theorem $2.5$} \ {\it For any point $p\in(M^n, {\mathcal
A}^{\omega})$ with a local chart $(U_p, \varphi_p)$,
$\varphi_p(p)=(x_1^,x_2^0,\cdots,x_n^0)$, if there are just $s$
euclidean directions along ${\bf e}_{i_1},{\bf
e}_{i_2},\cdots,{\bf e}_{i_s}$ for a point , then the dimension of
$T_p^*M^n$ is}

$${\rm dim}T_p^*M^n \ = \ 2n-s$$

\no{\it with a basis}

$$\{d x^{i_j}|_p \ | \ 1\leq j\leq s\}
\bigcup\{d^-x^l_p, d^+x^l|_p \ | \ 1\leq l\leq n \ and \
l\not=i_j, 1\leq j\leq s \},$$

\no{\it where}

$$dx^i|_p(\frac{\partial}{\partial x^j}|_p)=\delta_j^i \ and \
d^{\epsilon_i}x^i|_p(\frac{\partial^{\epsilon_i}}{\partial
x^j}|_p)=\delta_j^i
$$

\no{\it for $\epsilon_i\in\{+,-\}, 1\leq i\leq n$.}

\vskip 8mm

\no{\bf \S $3.$ Pseudo-Manifold Geometries}

\vskip 4mm

\no Here we introduce {\it Minkowski norms} on these
pseudo-manifolds $(M^n, {\mathcal A}^{\omega})$.

\vskip 4mm

\no{\bf Definition $3.1$} \ {\it A Minkowski norm on a vector
space $V$ is a function $F:V\rightarrow {\bf R}$ such that}

$(1)$ \ {\it $F$ is smooth on $V\backslash\{0\}$ and $F(v)\geq 0$
for $\forall v\in V$;}

$(2)$ \ {\it $F$ is $1$-homogenous, i.e., $F(\lambda v)=\lambda
F(v)$ for $\forall \lambda> 0$;}

$(3)$ \ {\it for all $y\in V\backslash\{0\}$, the symmetric
bilinear form $g_y: V\times V\rightarrow{\bf R}$ with}

$$g_y(u,v)=\sum\limits_{i,j}\frac{\partial^2 F(y)}{\partial y^i\partial y^j}$$

\no{\it is positive definite for $u,v\in V$.}

\vskip 3mm

Denote by $TM^n=\bigcup\limits_{p\in(M^n, {\mathcal
A}^{\omega})}T_pM^n$.

\vskip 4mm

\no{\bf Definition $3.2$} \ {\it A pseudo-manifold geometry is a
pseudo-manifold $(M^n, {\mathcal A}^{\omega})$ endowed with a
Minkowski norm $F$ on $TM^n$.}

\vskip 3mm

Then we get the following result.

\vskip 4mm

\no{\bf Theorem $3.1$} \ {\it There are pseudo-manifold
geometries.}

\vskip 3mm

{\it Proof} \ Consider an eucildean $2n$-dimensional space ${\bf
R}^{2n}$. Then there exists a Minkowski norm
$F(\overline{x})=|\overline{x}|$ at least. According to Theorem
$2.4$, $T_pM^n$ is ${\bf R}^{s+2(n-s)}$ if $\omega(p)$ has $s$
euclidean directions along ${\bf e}_1,{\bf e}_2,\cdots,{\bf e}_n$.
Whence there are Minkowski norms on each chart of a point in
$(M^n, {\mathcal A}^{\omega})$.

Since $(M^n,{\mathcal A})$ has  finite cover
$\{(U_{\alpha},\varphi_{\alpha})|\alpha\in I\}$, where $I$ is a
finite index set, by the decomposition theorem for unit, we know
that there are smooth functions $h_{\alpha}, \alpha\in I$ such
that

$$\sum\limits_{\alpha\in I}h_{\alpha}=1 \ {\rm with} \ 0\leq h_{\alpha}\leq 1.$$

Choose a Minkowski norm $F^{\alpha}$ on each chart
$(U_{\alpha},\varphi_{\alpha})$. Define

\[
F_{\alpha}=\left\{\begin{array}{cc}
h^{\alpha}F^{\alpha},& {\rm if}\quad p\in U_{\alpha} ,\\
0,& {\rm if}\quad p\not\in U_{\alpha}
\end{array}
\right.
\]

\no for $\forall p\in(M^n, \varphi^{\omega})$. Now let

$$F=\sum\limits_{\alpha\in I}F_{\alpha}.$$

\no Then $F$ is a Minkowski norm on $TM^n$ since it satisfies all
of these conditions $(1)-(3)$ in Definition $3.1$. \ \ $\natural$

Although the dimension of each tangent vector space maybe
different, we can also introduce {\it principal fiber bundles} and
{\it connections} on pseudo-manifolds.

\vskip 4mm

\no{\bf Definition $3.3$} \ {\it A principal fiber bundle (PFB)
consists of a pseudo-manifold $(P,{\mathcal A}_1^{\omega})$, a
projection $\pi:(P,{\mathcal A}_1^{\omega})\rightarrow(M,{\mathcal
A}_0^{\pi(\omega)})$, a base pseudo-manifold $(M,{\mathcal
A}_0^{\pi(\omega)})$ and a Lie group $G$, denoted by
$(P,M,\omega^{\pi},G)$ such that (1), (2) and (3) following hold.}

($1$) \ {\it There is a right freely action of $G$ on
$(P,{\mathcal A}_1^{\omega})$, i.e., for $\forall g\in G$, there
is a diffeomorphism $R_g:(P,{\mathcal
A}_1^{\omega})\rightarrow(P,{\mathcal A}_1^{\omega})$ with
$R_g(p^{\omega})=p^{\omega}g$ for $\forall p\in (P,{\mathcal
A}_1^{\omega})$ such that $p^{\omega}(g_1g_2)=(p^{\omega}g_1)g_2$
for $\forall p\in (P,{\mathcal A}_1^{\omega})$, $\forall
g_1,g_2\in G$ and $p^{\omega}e=p^{\omega}$ for some
$p\in(P^n,{\mathcal A}_1^{\omega})$, $e\in G$ if and only if $e$
is the identity element of $G$.}

($2$) \ {\it The map $\pi:(P,{\mathcal
A}_1^{\omega})\rightarrow(M,{\mathcal A}_0^{\pi(\omega)})$ is onto
with $\pi^{-1}(\pi(p))=\{pg| g\in G\}$, $\pi\omega_1=\omega_0\pi$,
and regular on spatial directions of $p$, i.e., if the spatial
directions of $p$ are $(\omega_1,\omega_2,\cdots,\omega_n)$, then
 $\omega_i$ and $\pi(\omega_i)$
are both elliptic, or euclidean, or hyperbolic and
$|\pi^{-1}(\pi(\omega_i))|$ is a constant number independent of
$p$ for any integer $i, 1\leq i\leq n$.}

($3$) \ {\it For $\forall x\in(M,{\mathcal A}_0^{\pi(\omega)})$
there is an open set $U$ with $x\in U$ and a diffeomorphism
$T_u^{\pi(\omega)}:(\pi)^{-1}(U^{\pi(\omega)})\rightarrow
U^{\pi(\omega)}\times G$ of the form
$T_u(p)=(\pi(p^{\omega}),s_u(p^{\omega}))$, where
$s_u:\pi^{-1}(U^{\pi(\omega)})\rightarrow G$ has the property
$s_u(p^{\omega}g)=s_u(p^{\omega})g$ for $\forall g\in G,
p\in\pi^{-1}(U)$.}

\vskip 3mm

We know the following result for principal fiber bundles of
pseudo-manifolds.

\vskip 4mm

\no{\bf Theorem $3.2$} \ {\it Let $(P,M,\omega^{\pi},G)$ be a PFB.
Then}

$$(P,M,\omega^{\pi},G)=(P,M,\pi,G)$$

\no{\it if and only if all points in pseudo-manifolds
$(P,{\mathcal A}_1^{\omega})$ are euclidean.}

\vskip 3mm

{\it Proof} \ For $\forall p\in(P,{\mathcal A}_1^{\omega})$, let
$(U_p,\varphi_p)$ be a chart at $p$. Notice that
$\omega^{\pi}=\pi$ if and only if $\varphi_p^{\omega}=\varphi_p$
for $\forall p\in(P,{\mathcal A}_1^{\omega})$. According to
Theorem $2.1$, by definition this is equivalent to that all points
in $(P,{\mathcal A}_1^{\omega})$ are euclidean. \ \ $\natural$

\vskip 4mm

\no{\bf Definition $3.4$} \ {\it Let $(P,M,\omega^{\pi},G)$ be a
PFB with ${\rm dim} G=r$. A subspace family
$H=\{H_p|p\in(P,{\mathcal A}_1^{\omega}), {\rm dim}H_p={\rm
dim}T_{\pi(p)}M \}$ of $TP$ is called a connection if conditions
($1$) and ($2$) following hold.}

($1$) \ {\it For $\forall p\in (P,{\mathcal A}_1^{\omega})$, there
is a decomposition}

$$T_pP=H_p\bigoplus V_p$$

\no{\it and the restriction $\pi_*|_{H_p}:H_p\rightarrow
T_{\pi(p)}M$ is a linear isomorphism.}

($2$) \ {\it $H$ is invariant under the right action of $G$, i.e.,
for $p\in (P,{\mathcal A}_1^{\omega})$, $\forall g\in G$,}

$$(R_g)_{*p}(H_p)=H_{pg}.$$

\vskip 3mm

Similar to Theorem $3.2$, the conception of connection introduced
in Definition $3.4$ is more general than the popular connection on
principal fiber bundles.

\vskip 4mm

\no{\bf Theorem $3.3$}(dimensional formula) \ {\it Let
$(P,M,\omega^{\pi},G)$ be a $PFB$ with a connection $H$. For
$\forall p\in (P,{\mathcal A}_1^{\omega})$, if the number of
euclidean directions of $p$ is $\lambda_P(p)$, then}

$${\rm dim}V_p=\frac{({\rm dim}P-{\rm dim}M)(2{\rm dim}P-\lambda_P(p))}{{\rm dim}P}.$$

\vskip 3mm

{\it Proof} \ Assume these euclidean directions of the point $p$
being ${\bf e}_1,{\bf e}_2,\cdots,{\bf e}_{\lambda_P(p)}$. By
definition $\pi$ is regular, we know that $\pi({\bf e}_1),
\pi({\bf e}_2),\cdots, \pi({\bf e}_{\lambda_P(p)})$ are also
euclidean in $(M,{\mathcal A}_1^{\pi(\omega)})$. Now since

$$\pi^{-1}(\pi({\bf e}_1))=\pi^{-1}(\pi({\bf e}_2))=\cdots =
\pi^{-1}(\pi({\bf e}_{\lambda_P(p)}))=\mu= \ {\rm constant},$$

\no we get that $\lambda_P(p)=\mu \lambda_M$, where $\lambda_M$
denotes the correspondent euclidean directions in $(M,{\mathcal
A}_1^{\pi(\omega)})$. Similarly, consider all directions of the
point $p$, we also get that ${\rm dim}P=\mu {\rm dim}M$.
Thereafter

$$\lambda_M = \frac{{\rm dim}M}{{\rm dim}P}\lambda_P(p). \ \ \ \ \ (3.1)$$

\no Now by Definition $3.4$, $T_pP=H_p\bigoplus V_p$, i.e.,

$${\rm dim}T_pP = {\rm dim}H_p+{\rm dim}V_p. \  \ (3.2)$$

\no Since $\pi_*|_{H_p}:H_p\rightarrow T_{\pi(p)}M$ is a linear
isomorphism, we know that ${\rm dim}H_p={\rm dim}T_{\pi(p)}M$.
According to Theorem $2.4$, we have formulae

$${\rm dim}T_pP= 2{\rm dim}P-\lambda_P(p)$$

\no and

$${\rm dim}T_{\pi(p)}M=2{\rm dim}M-\lambda_M=2{\rm dim}M-\frac{{\rm dim}M}{{\rm dim}P}\lambda_P(p).$$

\no Now replacing all these formulae into $(3.2)$, we get that

$$2{\rm dim}P-\lambda_P(p)=2{\rm dim}M-\frac{{\rm dim}M}{{\rm dim}P}\lambda_P(p)+{\rm dim}V_p.$$

\no That is,

$${\rm dim}V_p=\frac{({\rm dim}P-{\rm dim}M)(2{\rm dim}P-\lambda_P(p))}{{\rm dim}P}. \ \ \natural$$

We immediately get the following consequence by Theorem $3.3$.

\vskip 4mm

\no{\bf Corollary $3.1$} \ {\it Let $(P,M,\omega^{\pi},G)$ be a
$PFB$ with a connection $H$. Then for $\forall p\in (P,{\mathcal
A}_1^{\omega})$,}

$${\rm dim}V_p={\rm dim}P-{\rm dim}M$$

\no{\it if and only if the point $p$ is euclidean.}

\vskip 3mm

Now we consider conclusions included in Smarandache geometries,
particularly in pseudo-manifold geometries.

\vskip 4mm

\no{\bf Theorem $3.4$} \ {\it A pseudo-manifold geometry $(M^n,
\varphi^{\omega})$ with a Minkowski norm on $TM^n$ is a Finsler
geometry if and only if all points of $(M^n, \varphi^{\omega})$
are euclidean.}

\vskip 3mm

{\it Proof} \ According to Theorem $2.1$,
$\varphi_p^{\omega}=\varphi_p$ for $\forall p\in(M^n,
\varphi^{\omega})$ if and only if $p$ is eucildean. Whence, by
definition $(M^n, \varphi^{\omega})$ is a Finsler geometry if and
only if all points of $(M^n, \varphi^{\omega})$ are euclidean. \ \
$\natural$

\vskip 4mm

\no{\bf Corollary $3.1$} \ {\it There are inclusions among
Smarandache geometries, Finsler geometry, Riemann geometry and
Weyl geometry:}

\begin{eqnarray*}
& \ &\{Smarandache \ geometries\}\supset\{pseudo-manifold \
geometries\}\\
& \ &\supset\{Finsler \ geometry\}\supset\{Riemann \
geometry\}\supset\{Weyl \ geometry\}.
\end{eqnarray*}

\vskip 3mm

{\it Proof} \ The first and second inclusions are implied in
Theorems $2.1$ and $3.3$. Other inclusions are known in a
textbook, such as $[4]-[5]$. \ \ $\natural$

Now we consider complex manifolds. Let $z^i=x^i+\sqrt{-1}y^i$. In
fact, any complex manifold $M^{n}_c$ is equal to a smooth real
manifold $M^{2n}$ with a natural base $\{\frac{\partial}{\partial
x^i}, \frac{\partial}{\partial y^i}\}$ for $T_pM^n_c$ at each
point $p\in M^n_c$. Define a {\it Hermite} manifold $M^{n}_c$ to
be a manifold $M^{n}_c$ endowed with a Hermite inner product
$h(p)$ on the tangent space $(T_pM^n_c, J)$ for $\forall p\in
M^n_c$, where $J$ is a mapping defined by

$$J(\frac{\partial}{\partial x^i}|_p)=\frac{\partial}{\partial y^i}|_p,
\ \ J(\frac{\partial}{\partial y^i}|_p)=-\frac{\partial}{\partial
x^i}|_p
$$

\no at each point $p\in M^n_c$ for any integer $i, 1\leq i\leq n$.
Now let

$$h(p)=g(p)+\sqrt{-1}\kappa(p), \ \ p\in M^m_c.$$

\no Then a {\it K\"{a}hler manifold} is defined to be a Hermite
manifold $(M^n_c,h)$ with a closed $\kappa$ satisfying

$$\kappa(X,Y)=g(X,JY), \ \forall X,Y\in T_pM^n_c, \forall p\in M^n_c.$$

Similar to Theorem $3.3$ for real manifolds, we know the next
result.

\vskip 4mm

\no{\bf Theorem $3.5$} \ {\it A pseudo-manifold geometry $(M^n_c,
\varphi^{\omega})$ with a Minkowski norm on $TM^n$ is a K\"{a}hler
geometry if and only if $F$ is a Hermite inner product on $M^n_c$
with all points of $(M^n, \varphi^{\omega})$ being euclidean.}

\vskip 3mm

{\it Proof} \ Notice that a complex manifold $M^n_c$ is equal to a
real manifold $M^{2n}$. Similar to the proof of Theorem $3.3$, we
get the claim. \ \ $\natural$

As a immediately consequence, we get the following inclusions in
Smarandache geometries.

\vskip 4mm

\no{\bf Corollary $3.2$}\ {\it There are inclusions among
Smarandache geometries, pseudo-manifold geometry and K\"{a}hler
geometry:}

\begin{eqnarray*}
\{Smarandache \ geometries\}&\supset&\{pseudo-manifold \
geometries\}\\ &\supset&\{K\ddot{a}hler \ geometry\}.
\end{eqnarray*}

\vskip 8mm

\no{\bf \S $4.$ Further Discussions}

\vskip 4mm

\no Undoubtedly, there are many and many open problems and
research trends in pseudo-manifold geometries. Further research
these new trends and solving these open problems will enrich one's
knowledge in sciences.

Firstly, we need to get these counterpart in pseudo-manifold
geometries for some important results in Finsler geometry or
Riemann geometry.

\vskip 3mm

\no{\bf $4.1.$ Storkes Theorem} \ {\it Let $(M^n,{\mathcal A})$ be
a smoothly oriented manifold with the $T_2$ axiom hold. Then for
$\forall\varpi\in A_0^{n-1}(M^n)$,}

$$\int_{M^n}d\varpi=\int_{\partial M^n}\varpi.$$

\no This is the well-known Storkes formula in Riemann geometry. If
we replace $(M^n,{\mathcal A})$ by $(M^n,{\mathcal A}^{\omega})$,
{\it what will happens?} Answer this question needs to solve
problems following.\vskip 3mm

($1$) \ {\it Establish an integral theory on pseudo-manifolds.}

($2$) \ {\it Find conditions such that the Storkes formula hold
for pseudo-manifolds.}

\vskip 3mm

\no{\bf $4.2.$ Gauss-Bonnet Theorem} \ {\it Let $S$ be an
orientable compact surface. Then}

$$\int\int_SKd\sigma= 2\pi\chi(S),$$

\no {\it where $K$ and $\chi(S)$ are the Gauss curvature and Euler
characteristic of $S$} This formula is the well-known Gauss-Bonnet
formula in differential geometry on surfaces. Then {\it what is
its counterpart in pseudo-manifold geometries?} This need us to
solve problems following.\vskip 3mm

($1$) \ {\it Find a suitable definition for curvatures in
pseudo-manifold geometries.}

($2$) \ {\it Find generalizations of the Gauss-Bonnect formula for
pseudo-manifold geometries, particularly, for
pseudo-surfaces.}\vskip 2mm

For a oriently compact Riemann manifold $(M^{2p},g)$, let

$$\Omega=\frac{(-1)^p}{2^{2p}\pi^pp!}\sum\limits_{i_1,i_2,\cdots,i_{2p}}
\delta_{1,\cdots,2p}^{i_1,\cdots,i_{2p}}\Omega_{i_1i_2}\wedge\cdots\wedge
\Omega_{i_{2p-1}i_{2p}},$$

\no where $\Omega_{ij}$ is the curvature form under the natural
chart $\{e_i\}$ of $M^{2p}$ and

\[
\delta_{1,\cdots,2p}^{i_1,\cdots,i_{2p}}=\left\{\begin{array}{cc}
1,& {\rm if \ permutation} \ i_1\cdots i_{2p} \ {\rm is \ even,}\\
-1,& {\rm if \ permutation} \ i_1\cdots i_{2p} \ {\rm is \ odd,}\\
0, & {\rm otherwise.}
\end{array}
\right.
\]

Chern proved that$^{[4]-[5]}$

$$\int_{M^{2p}}\Omega = \chi(M^{2p}).$$

\no Certainly, these new kind of global formulae for
pseudo-manifold geometries are valuable to find.

\vskip 3mm

\no{\bf $4.3.$ Gauge Fields} \ Physicists have established a gauge
theory on principal fiber bundles of Riemann manifolds, which can
be used to unite gauge fields with gravitation. Similar
consideration for pseudo-manifold geometries will induce new gauge
theory, which enables us to asking problems following.

{\it Establish a gauge theory on those of pseudo-manifold
geometries with some additional conditions.}

($1$) \ {\it Find these conditions such that we can establish a
gauge theory on a pseudo-manifold geometry.}

($2$) \ {\it Find the Yang-Mills equation in a gauge theory on a
pseudo-manifold geometry.}

($2$) \ {\it Unify these gauge fields and gravitation.}

\vskip 10mm

\no{\bf References}\vskip 5mm

\re{[1]}R.Abraham, J.E.Marsden and T.Ratiu, {\it Manifolds, tensor
analysis, and applications}, Addison-Wesley Publishing Company,
Inc. 1983.

\re{[2]}G.Bassini and S.Capozziello, Multi-Spaces and many worlds
from conservation laws, {\it Progress in Physics}, Vol.4(2006),
65-72.

\re{[3]}D.Bleecker, {\it Gauge theory and variational principles},
Addison-Wesley Publishing Company, Inc. 1981.

\re{[4]}S.S.Chern and W.H.Chern, {\it Lectures in Differential
Geometry}(in Chinese), Peking University Press, 2001.

\re{[5]}W.H.Chern and X.X.Li, {\it Introduction to Riemann
Geometry}, Peking University Press, 2002.

\re{[6]}H.Iseri, {\it Smarandache manifolds}, American Research
Press, Rehoboth, NM,2002.

\re{[7]}H.Iseri, {\it Partially Paradoxist Smarandache
Geometries}, http://www.gallup.unm.
edu/\~smarandache/Howard-Iseri-paper.htm.

\re{[8]}L.Kuciuk and M.Antholy, An Introduction to Smarandache
Geometries, {\it Mathematics Magazine, Aurora, Canada},
Vol.12(2003).

\re{[9]}L.F.Mao, {\it On Automorphisms groups of Maps, Surfaces
and Smarandache geometries}, {\it Sientia Magna}, Vol.$1$(2005),
No.$2$, 55-73.

\re{[10]}L.F.Mao, A new view of combinatorial maps by
Smarandache's notion, {\it arXiv: math.GM/0506232}.

\re{[11]}L.F.Mao, {\it Automorphism Groups of Maps, Surfaces and
Smarandache Geometries}, American Research Press, 2005.

\re{[12]}L.F.Mao, {\it Smarandache multi-space theory}, Hexis,
Phoenix, AZ£¬2006.

\re{[13]}L.F.Mao, On multi-metric spaces, {\it Scientia Magna},
Vol.2, No.1(2006), 87-94.

\re{[14]}F.Smarandache, {\it A Unifying Field in Logics.
Neutrosopy: Neturosophic Probability, Set, and Logic}, American
research Press, Rehoboth, 1999.

\re{[15]}F.Smarandache, A Unifying Field in Logic: Neutrosophic
Field, {\it Multi-Valued Logic}, Vol.8, No.3(2002)(special issue
on Neutrosophy and Neutrosophic Logic), 385-438.

\re{[16]}F.Smarandache, Mixed noneuclidean geometries, {\it eprint
arXiv: math/0010119}, 10/2000.

\end{document}